\documentclass[12pt,a4paper]{article}
\usepackage{amsmath}
\usepackage{amsfonts}
\usepackage{amssymb}
\usepackage{graphicx}
\newtheorem{theorem}{Theorem}
\newtheorem{lemma}{Lemma}
\newtheorem{proposition}{Proposition}
\newtheorem{remark}{Remark}

\newenvironment{proof}[1][Proof]{\noindent\textbf{#1.} }{\ \rule{0.5em}{0.5em}}

\begin{document}

\title{Similarity solutions for high frequency excitation of liquid metal in an antisymmetric magnetic field}
\date{}
\author{}
\maketitle
\vspace{-1cm}

\begin{center}Bernard BRIGHI\end{center}

\vspace{-.7cm}
\begin{center}Laboratoire de Math\'ematiques, Informatique et Applications, Universit\'e de Haute-Alsace\\
4 rue des fr\`eres Lumi\`ere, 68093 Mulhouse, FRANCE\\
E-mail:  bernard.brighi@uha.fr
\end{center}

\vspace{.1cm}

\begin{center}Jean-David HOERNEL\end{center}

\vspace{-.7cm}
\begin{center}Department of Mathematics, Technion-Israel Institute of Technology\\
Amado Bld., Haifa, 32000 ISRAEL\\
E-mail: j-d.hoernel@wanadoo.fr
\end{center}

\footnotetext{AMS 2000 Subject Classification: 34B15, 34C11, 76D10.}
\footnotetext{Key words and phrases: Third order differential equations, boundary value problems, plane dynamical systems, blowing-up co-ordinates.}
\footnotetext{The second author thanks the Department of Mathematics of the Technion for supporting his researchs through a Postdoctoral Fellowship.}

\abstract{The aim of this paper is to investigate, as precisely as possible, a boundary value problem involving a third order ordinary differential equation. Its solutions are the similarity solutions of a problem arising in the study of the phenomenon of high frequency excitation of liquid metal systems in an antisymmetric magnetic field within the framework  of boundary layer approximation.}

\section{Introduction}
In this paper we study the third order non-linear autonomous differential
equation
\begin{equation}
f'''+\frac{m+1}{2}ff''-m f'^2=0 \label{equation}
\end{equation}
on $[0,\infty)$, with the boundary conditions
\begin{align}
f(0)&=a, \label{cond01}\\
f'(0)&=-1, \label{cond02} \\
f'(\infty)& =0 \label{cond03}
\end{align}
where $m\in\mathbb{R}$, $a\in\mathbb{R}$ and $f'(\infty):=\underset{t\rightarrow\infty}{\lim} f'(t)$. 

This boundary value problem appears in the paper of H. K. Moffatt \cite{mof} and is related to the behavior of a liquid metal in an antisymmetric field, in the framework of boundary layer approximation.

The study of similarity solutions for free convection in a fluid saturated porous medium near a semi-infinite vertical flat plate on which the heat is precribed or high frequency excitation of liquid metal systems in an symmetric magnetic field, both in the framework of boundary layer approximation, leads to the same third order ordinary differential equation (\ref{equation}) subjected to the boundary conditions $f(0)=a$, $f'(0)=1$ and $f'(\infty)=0.$
This problem also appears when studying boundary layer flows adjacent to stretching walls. 
One can find explicit solutions of this problem for some particular values of $m$ in \cite{brighicr}, \cite{brighi02}, \cite{brighi01}, \cite{crane}, \cite{gup},  \cite{ing}, \cite{mag} and \cite{stu}. For mathematical results about existence, nonexistence, uniqueness, nonuniqueness and asymptotic behavior, see \cite{banks1}, \cite{brighicr}, \cite{brighi02} and \cite{ing} for $a=0$, and  \cite{brighi01},  \cite{equiv}, \cite{brighisari}, \cite{guedda} and  \cite{guedda1} for the general case. Numerical  investigations can be found in \cite{banks1}, \cite{brighi04}, \cite{pop1}, \cite{cheng}, \cite{ing}, \cite{mag} and \cite{wood}. For the high frequency excitation of liquid metal systems in an symmetric magnetic field, see \cite{mof}.

When studying similarity solutions for free convection in a fluid saturated porous medium near a semi-infinite vertical flat plate on which the heat flux is precribed again in the framework of boundary layer approximation, we obtain this time the equation $f'''+(m+2)ff''-(2m+1)f'^2=0$ which differs from (\ref{equation}) only by its coefficients, with the boundary conditions $f(0)=a$, $f''(0)=-1$ and $f'(\infty)=0.$ Numerical results can be found in \cite{pop} and the mathematical study of existence, uniqueness and qualitative properties of the solutions of this problem is made in \cite{heat_flux}. 

For a survey of the previously described problems, see \cite{gaeta}. 

One particular case of all these equations is the Blasius equation $f'''+ff''=0$ introduced in \cite{bla}.  The Blasius equation is obtained by setting $m=0$ and doing some proper rescaling in (\ref{equation}). The corresponding problem with the boundary conditions  $f(0)=a$, $f'(0)=b\geq 0$ and $f'(\infty)=\lambda$ admits an unique solution for $\lambda\geq 0$, and no solution for $\lambda<0$. This well known case is studied, for example, in 
\cite{brighi03},  \cite{coppel} and \cite{hart}. 
On the other hand, with the boundary conditions (\ref{cond01})-(\ref{cond03}), the situation is completely different. In fact, one can show that for  $a=\sqrt{3}$, the Blasius problem
$$\left\{
\begin{array}[c]{l}
f'''+ff''=0, \\ \noalign{\vskip2mm}
f(0)=a,\ f'(0)=-1,\  f'(\infty)=0
\end{array}
\right.$$
admits infinitely many solutions, and for every $n\in\mathbb{N}$, there are values of $a$ such that this problem has exactly $n$ solutions. See \cite{bfs} for the proofs of these results. In the remainder of the paper we will only consider $m\neq 0$. 

The study of similarity solutions for mixed convection in a fluid saturated porous medium near a semi-infinite vertical flat plate on which the heat is precribed, leads to the equation $f'''+(m+1)ff''+2m(1-f')f'=0$ with the boundary conditions $f(0)=a$, $f'(0)=b$ and $f'(\infty)=1$. Results about it can be found in \cite{aly}, \cite{aml}, \cite{guedda2} and \cite{nazar}. 

The Falkner-Skan equation $f'''+ff''+m(1-f'^2)=0$ is in the same family of problems. See, for example, \cite{coppel}, \cite{falk}, \cite{hart}, \cite{ish1}, \cite{wang}, \cite{yang} and \cite{yang2} for results about it.

New results about the more general equation $f'''+ff''+g(f')=0$ with the boundary conditions $f(0)=a$, $f'(0)=b$ and $f'(\infty)=c$ for some given function $g$ can be found in \cite{jde}, see also \cite{utz}.

For some new results about the full model of free convection in a plane and bounded fluid saturated porous medium, see \cite{ake}.

\section{Preliminary results} 
First of all, let us notice that for every $\tau>0$, the function $t\mapsto\frac{6}{t+\tau}$ is a solution of the equation (\ref{equation}) for any value of $m$, and thus for $a=\sqrt{6}$, the function $f(t)=\frac{6}{t+\sqrt{6}}$ is a solution of the problem (\ref{equation})-(\ref{cond03}).

Now, we remark that, if $f$ verifies (\ref{equation}), then
\begin{equation}
\left(f''e^{\frac{m+1}{2}F}\right)'=mf'^2e^{\frac{m+1}{2}F} \label{egalite01}
\end{equation}
with $F$ any anti-derivative of $f.$ 
\medskip
\begin{lemma}
\label{concavite}Let $f$ be a non constant solution of the equation $(\ref{equation})$ on some interval $I$. For all $t_{0}$ in $I$ we have that
\begin{itemize}
\item if $m<0$, $f''(t_{0})\leq0\Rightarrow f''(t)<0$ for $t>t_{0},$
\item if $m>0$, $f''(t_{0})\geq0\Rightarrow f''(t)>0$ for $t>t_{0}.$
\end{itemize}
\end{lemma}
\begin{proof} Immediate using (\ref{egalite01}) and the fact that  $f^{\prime}$ and $f^{\prime\prime}$
cannot vanish at the same point without being identically equal to zero.\end{proof}
\medskip

\begin{proposition}
\label{prop<0}Let $m<0$. If $f$ is a solution of the problem $(\ref{equation})$-$(\ref{cond03})$ then $f''(0)>0$, and
\begin{itemize}
\item either $f$ is convex and decreasing on $[0,\infty)$,
\item or there exists $t_0$ with $f''(t_0)=0$ and $f'(t_0)\geq 0$ such that $f$ is convex and first decreasing then increasing on $[0,t_0)$, and concave and increasing on $[t_0,\infty)$.  Moreover, $f$ is negative at infinity for $m\leq-1$, and
positive at infinity for $-1\leq m<0$ $($and in particular this implies that such solutions cannot exist for $m=-1)$.
\end{itemize}
\end{proposition}
\begin{proof}
Suppose that $f''(0)\leq 0$, then, using Lemma \ref{concavite}, we have $f''<0$ and $f'$ is decreasing. This is a contradiction with $f'(0)=-1$ and $f'(\infty)=0$.
Hence $f''(0)>0$. 

If $f''$ never vanishes, then $f$ is convex and $f'$ is increasing. Since $f'(0)=-1$ and $f'(\infty)=0$ we get $-1\leq f'<0$. 

If there exists a $t_0$ such that $f''>0$ on $[0,t_0)$ and $f''(t_0)=0$ then, by Lemma \ref{concavite} we have $f''<0$ on $(t_0,\infty)$. Hence $f'$ is decreasing on $(t_0,\infty)$ and since $f'(\infty)=0$, we must have $f'>0$ on $[t_0,\infty)$ and $f$ is increasing on $[t_0,\infty)$. 

For $-1\leq m <0$, if $f$ is negative at infinity, there exists $t_1$ such that $f<0$, $f'>0$ and $f''<0$ on $(t_1,\infty)$. Therefore
$$f'''=mf'^2-\frac{m+1}{2}ff''< -\frac{m+1}{2}ff''\leq 0$$
on $(t_1,\infty)$.
This implies that $f'$ is concave on $(t_1,\infty)$, a contradiction with the facts that $f'>0$ on $(t_1,\infty)$ and $f'(\infty)=0$. For $m\leq -1$, the same arguments shows that $f$ is negative at infinity.\end{proof}
\vskip4pt plus2pt

\begin{remark} In  {\rm \cite{guedda1}} $(${\rm Theorem 2.1}$)$, it is proved that if $-1<m<0$, then any solution $f$ of $(\ref{equation})$-$(\ref{cond02})$ such that $f''(0)<0$ only exists on $[0,T)$ with $0<T<\infty$ and that $\underset{t\rightarrow T}{\lim} f(t)=-\infty$.
\end{remark}

\medskip

\begin{proposition}
\label{prop>0}Let $m> 0$. If $a\leq 0$, there are no solutions of the problem $(\ref{equation})$-$(\ref{cond03})$. If $a>0$, and if $f$ is a solution of $(\ref{equation})$-$(\ref{cond03})$ then
$f$ is positive, decreasing and moreover
\begin{itemize}
\item if $f''(0)\geq 0$, then $f$ is convex,
\item if $f''(0)<0$, then there exists $t_0>0$ such that $f$ is concave on $[0,t_0]$ and convex on $(t_0,\infty)$.
\end{itemize}
\end{proposition}
\begin{proof} Let $f$ be a solution of $(\ref{equation})$-$(\ref{cond03})$.
First, let us suppose that $f''(0)\geq 0$. By Lemma \ref{concavite} we have that $f''>0$ everywhere. Hence  $f$ is convex and  $f'$ is increasing. Since $f'(0)=-1$ and $f'(\infty)=0$, we get $-1\leq f'\leq 0$ and $f$ is decreasing on $[0,\infty)$.

Now, let us suppose that $f''(0)<0$. If $f''<0$ on $[0,\infty)$, then $f'$ is decreasing and as $f'(0)=-1$ we cannot have $f'(\infty)=0$. Thus there exists $t_0>0$ such that $f''<0$ on $[0,t_0)$ and $f''(t_0)=0$.  
Then, $f'$ is decreasing on $[0,t_0]$ and we have $f'(t_0)\leq f'\leq -1$ on $[0,t_0]$. 
Moreover, by Lemma \ref{concavite} we get $f''>0$ on $(t_0,\infty)$, and $f'$ is increasing on $(t_0,\infty)$. As $f'(t_0)\leq -1$  we get $f'(t_0)\leq f'<0$ on $[t_0,\infty)$, therefore $f$ is decreasing on $[0,\infty)$. 

As $m>0$, if $f<0$ at infinity, then
$$f'''= -\frac{m+1}{2}ff''+mf'^2\geq mf'^2\geq 0$$
and $f'$ is convex at infinity. But as $f'<0$ and $f'(\infty)=0$ this cannot be the case.
Hence, $f>0$ at infinity, and since $f$ is decreasing, we get $f>0$ on $[0,\infty)$. This, in particular, implies that $a>0$, and the proof is complete.\end{proof}
\medskip

\section{Useful tools} In this part, we first give some identities and properties about solutions of $(\ref{equation})$, and next introduce blowing-up co-ordinates associated to $(\ref{equation})$ and related to the fact that if $f$ is a solution of $(\ref{equation})$, then it is also the case for the function $t\mapsto\kappa f(\kappa t)$.

\medskip
Let $f$ be a solution of $(\ref{equation})$ on some interval $[\alpha,\beta]$.
Integrating the equation $(\ref{equation})$ between $\alpha$ and $\beta$ leads to
\begin{equation}
f''(\beta)-f''(\alpha)+\frac{m+1}{2}f(\beta)f'(\beta)-\frac{m+1}{2}f(\alpha)f'(\alpha)=\frac{3m+1}{2}\int_\alpha^\beta f'^2(t)dt. \label{i1}
\end{equation}
Multiplying the equation $(\ref{equation})$ by $f$ and integrating between $\alpha$ and $\beta$ leads to
\begin{multline}
f(\beta)f''(\beta)-f(\alpha)f''(\alpha)-\frac{1}{2}f'^2(\beta)+\frac{1}{2}f'^2(\alpha)+
\frac{m+1}{2}f^2(\beta)f'(\beta)
\\
-\frac{m+1}{2}f^2(\alpha)f'(\alpha)=(2m+1)\int_\alpha^\beta f(t)f'^2(t)dt. \label{i3}
\end{multline}
Multiplying the equation $(\ref{equation})$ by $f''$ and integrating between $\alpha$ and $\beta$ leads to
\begin{equation}
\frac{1}{2}f''^2(\beta)-\frac{1}{2}f''^2(\alpha)-\frac{m}{3}f'^3(\beta)+\frac{m}{3}f'^3(\alpha)=-\frac{m+1}{2}\int_\alpha^\beta f(t)f''^2(t)dt. \label{i2}
\end{equation}
\medskip
\begin{proposition}\label{f2}
Let $m\in \mathbb R$. If $f$ is a solution of the problem $(\ref{equation})$-$(\ref{cond03})$ then we have
\begin{equation}
\underset{t\rightarrow \infty}{\lim} f''(t)=0\label{f''}
\end{equation}
and, if $m\neq -1$, there exists a sequence $t_n\rightarrow \infty$ such that
\begin{equation}
\underset{n\rightarrow \infty}{\lim} f'''(t_n)=\underset{n\rightarrow \infty}{\lim} f(t_n)f''(t_n)=0.\label{f'''}
\end{equation}
\end{proposition}
\begin{proof} Since $f'(\infty)=0$, there exists an increasing sequence $s_n$ such that
$s_n\to\infty$ and $f''(s_n)\to 0$ as $n\to\infty$. But, using $(\ref{i3})$, we see that $f''^2$ has a limit at infinity and hence $(\ref{f''})$ holds. In addition, choosing $t_n$ such that $f'''(t_n)=f''(n+1)-f''(n)$, and using  $(\ref{equation})$ we get $(\ref{f'''})$.\end{proof}

\medskip

Consider now a right maximal interval $I=\left[0,T\right)$ on which
$f$ does not vanish. For all $t$ in $I$, set
\begin{equation}
s=\int_{0}^{t}f(\xi)d\xi,\quad u(s)=\frac{f^{\prime}(t)}{f^2(t)}, \quad
v(s)=\frac{f^{\prime\prime}(t)}{f^3(t)}
\label{new_function}
\end{equation}
to obtain the system
\begin{equation}
\left\{
\begin{array}
[c]{l}
\dot{u}=P(u,v):=v-2u^{2},\\
\dot{v}=Q_{m}(u,v):=-\frac{m+1}{2}v+mu^{2}-3uv
\end{array}
\right.  \label{system}
\end{equation}
in which the dot denotes the differentiation with respect to $s.$ Let us notice that if $f$ is negative on $I$ then $s$ decreases as $t$ grows.

The singular points of $\left(  \ref{system}\right)  $ are $O=(0,0)$ and
$A=\left(  -\frac{1}{6},\frac{1}{18}\right)  .$ The isoclinic curves $P(u,v)=0$ and
$Q_{m}(u,v)=0$ are given by $v=2u^{2}$ and $v=\Psi_{m}(u)$ where
\[
\Psi_{m}(u)=\dfrac{mu^{2}}{3u+\frac{m+1}{2}}.
\]

The point $A$ is an unstable node for $m\leq\frac{4-2\sqrt{6}}{3}$, an unstable focus if $\frac{4-2\sqrt{6}}{3}<m<\frac{4}{3}$, 
a stable focus if $\frac{4}{3}<m<\frac{4+2\sqrt{6}}{3}$ and a stable node if $m\geq\frac{4+2\sqrt{6}}{3}$.

For $m\neq-1,$ the singular point
$O$ is a saddle-node of multiplicity 2. It admits a center manifold $\mathcal W_{0}$ that is
tangent to the subspace $L_{0}={\rm Sp}\left\{  (1,0)\right\}$ and a stable (resp. unstable) 
manifold $\mathcal W$ if $m>-1$ (resp $m<-1)$ that is tangent to the subspace 
$L={\rm Sp}\left\{  (1,-\frac{m+1}{2})\right\}.$ In the neighborhood of $O$, the manifold ${\mathcal W}$ takes place below $L$ when $m<-1$ or $m>-\frac{1}{3}$ and above $L$ when $-1<m<-\frac{1}{3}$. In the neighborhood of $O$, the center manifold ${\mathcal W}_0$ takes place above $L_0$ when $m<-1$ or $m>0$, and below $L_0$ when $-1<m<0$. 

We will not precise the behavior of the manifolds $\mathcal W$ and $\mathcal W_0$ for $m=-1$, because we will not use the co-ordinates $u$ and $v$ in this case.

In order to describe the phase portrait of the vector field in the neighborhood of the
saddle-node $O$ we will assume that the parabolic sector is delimited by the separatrices
$S_0$ and $S_1$ which are tangent to $L$, and the hyperbolic sectors are delimited, one by the separatrix $S_0$ and the separatrix $S_2$, which is tangent to $L_0$, and the other one by the separatrices $S_1$ and $S_2$. 
With these notations, we have that 
$${\mathcal W}=S_1 \cup \{O\} \cup S_0 \quad \text{and} \quad {\mathcal W}_0=S_2 \cup \{O\} \cup C_3$$
where $C_3$ is a phase curve.

We will use the superscript $+$ for $\omega$-separatrices and $-$ for 
$\alpha$-separatrices to obtain the behaviors
described  in the figure 1.

\medskip
\begin{center}
\begin{tabular}{cc}
\includegraphics[scale=0.28]{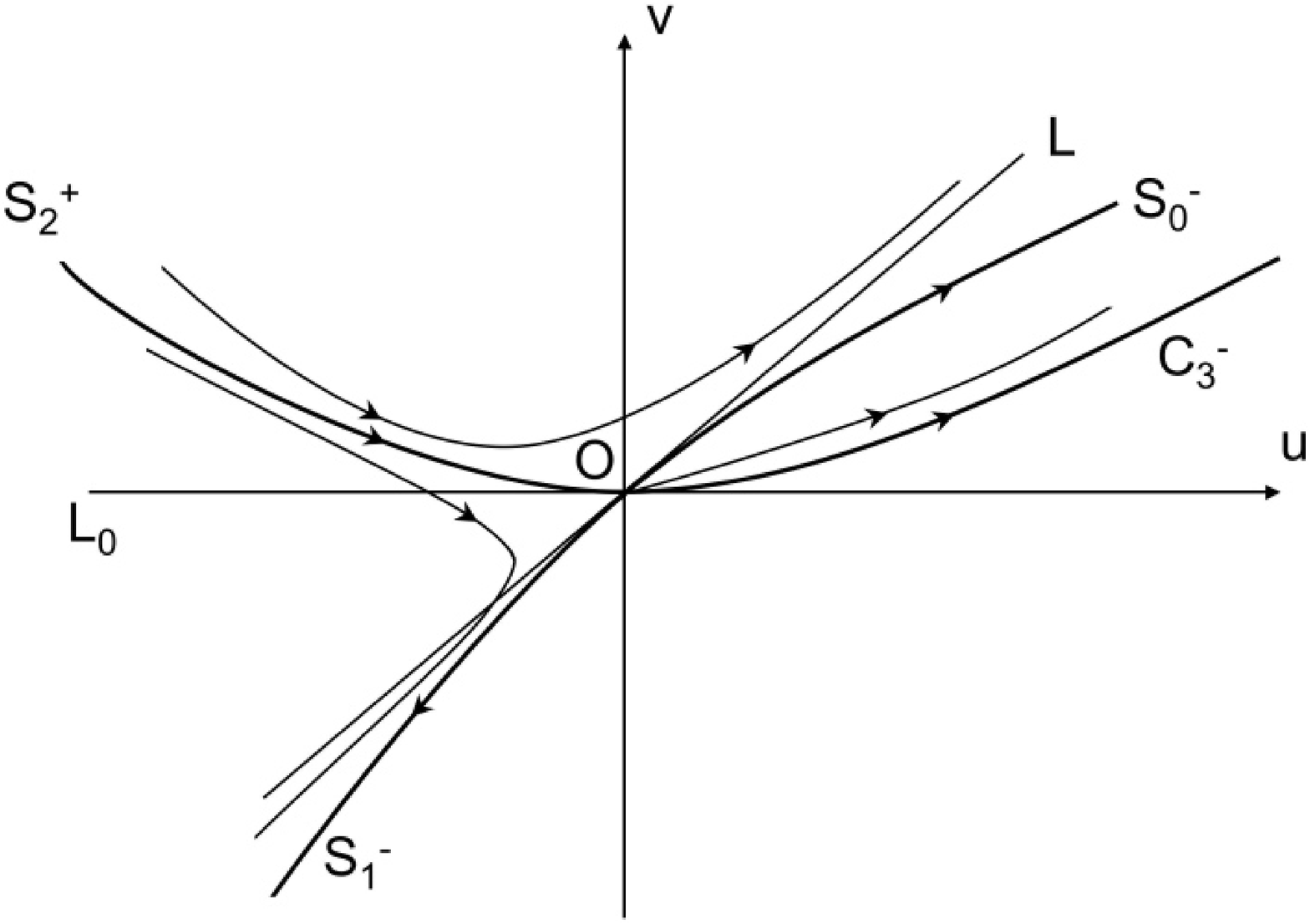}
&
\includegraphics[scale=0.28]{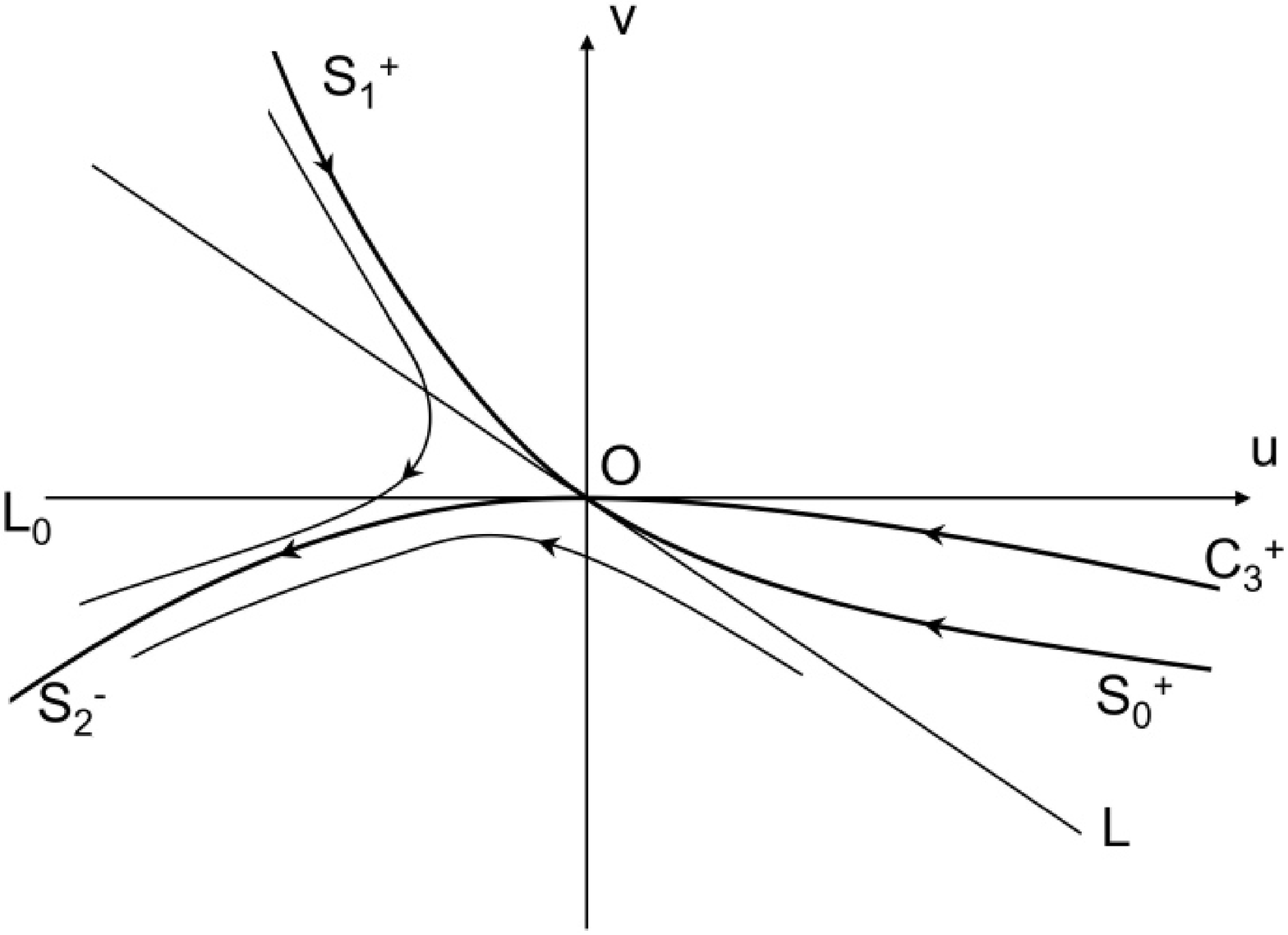}
\\
$m<-1$
&
$-1<m<-1/3$
\end{tabular}
\begin{tabular}{cc}
\includegraphics[scale=0.28]{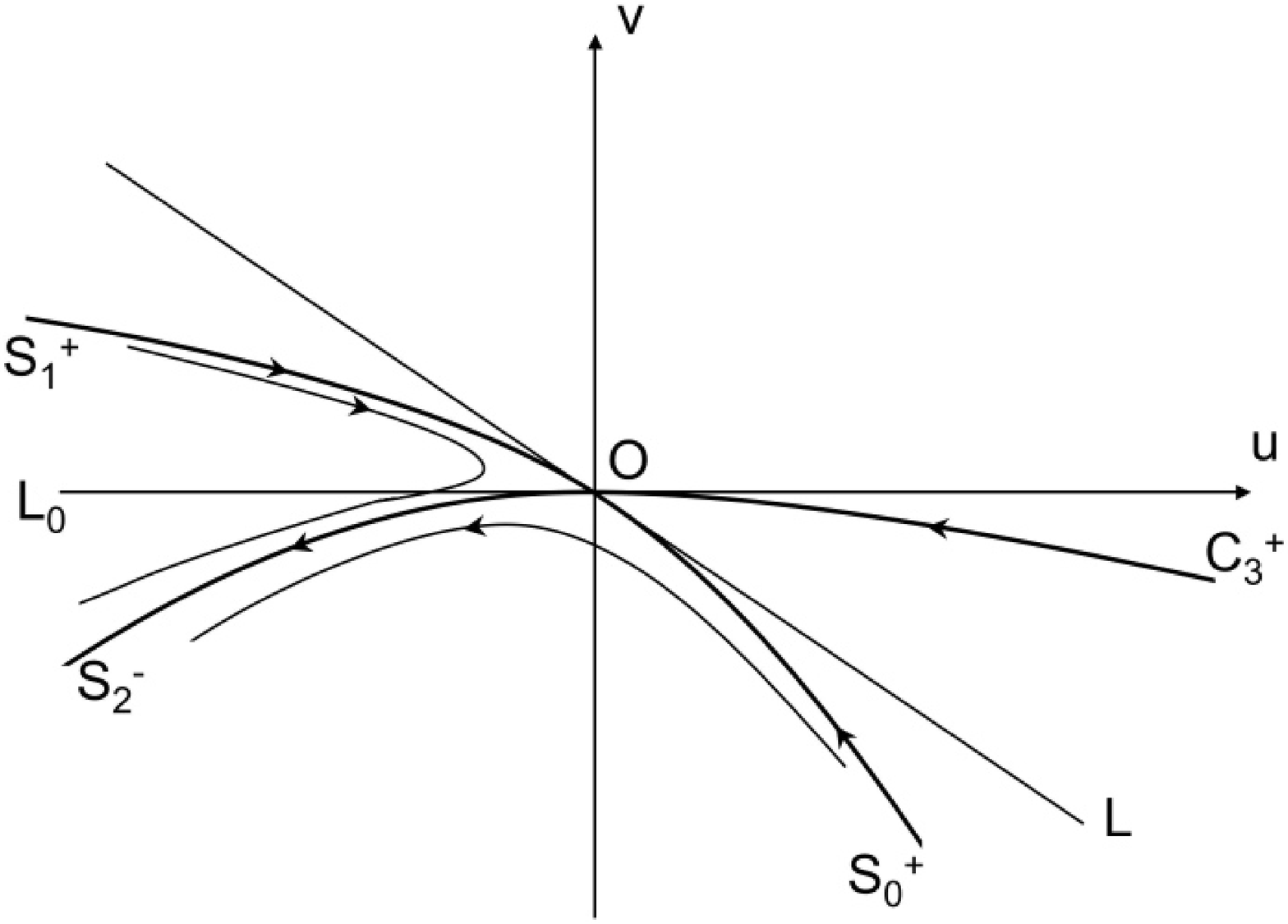}
&
\includegraphics[scale=0.28]{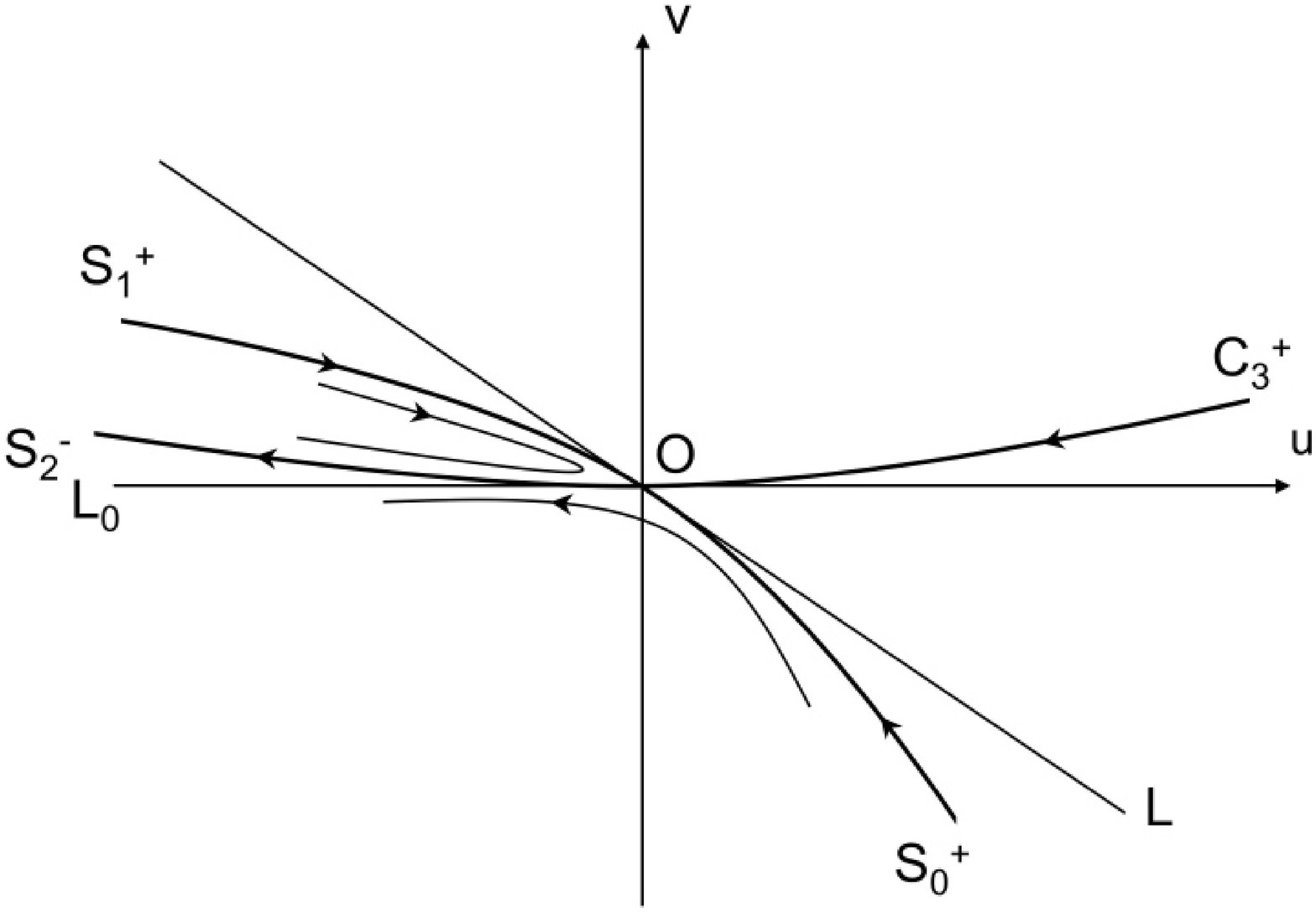}
\\
$-1/3<m<0$
&
$m>0$
\end{tabular}
\end{center}
\centerline{Figure 1}

\medskip
To study the global behavior of the separatrices, consider any connected piece of a phase curve $C$ of the plane dynamical system $(\ref{system})$ lying in the region
$P(u,v)<0$ (resp. $P(u,v)>0$); then $C$ can be characterized by $v=V_m(u)$
(resp. $v=W_m(u)$) with $u$ belonging to some interval, and  where $V_m$ (resp.
$W_m$) is a solution of the differential equation
\begin{equation}
\frac{dv}{du}=F_m(u,v):={Q_m(u,v)\over P(u,v)}={-{m+1\over 2}v+m u^2-3uv\over v-2u^2}.\label{eq9}
\end{equation}

To deduce results about the original problem $(\ref{equation})$-$(\ref{cond03})$, most of the time, we will consider the initial value problem
$$\left\{
\begin{array}[c]{l}
f'''+\frac{m+1}{2}ff''-mf'^2=0, \\
f(0)=a, \\
f'(0)=-1, \\
f''(0)=b
\end{array}
\right.
\leqno({\mathcal P}_{m,a,b})
$$
with $a\not=0$, and look at the trajectory $C_{a,b}$ for a given value of $m$ of the plane dynamical system $(\ref{system})$ defined by $(\ref{new_function})$ for 
$$u(0)=-\frac{1}{a^2}~~~~\hbox{ and }~~~~v(0)=\frac{b}{a^3}.$$
\newpage
\section{Main results}

\subsection{The case $m<-1$.}
\begin{theorem} \label{cv1}
Let $m<-1$. For every $a\in\Bbb R$, the problem $(\ref{equation})$-$(\ref{cond03})$ has a unique convex solution. Moreover, this solution is bounded, and if $\underset{t\rightarrow \infty}{\lim} f(t)=\ell$, then we have:  $-\sqrt{a^2-\frac4{m+1}}<\ell<a$.
\end{theorem} 
\begin{proof} For $f:[0,\infty)\to\Bbb R$ let us set $\tilde f:[0,\infty)\to\Bbb R$ defined by
$$\tilde f(s)=-\textstyle\sqrt{-\textstyle\frac{m+1}2}f\left(\sqrt{-\frac2{m+1}}~\!s\right).$$
Easily, one see that $f$ is a convex solution of  $(\ref{equation})$-$(\ref{cond03})$ if and only if $\tilde f$ is a concave solution of the problem
\begin{equation}
\left\{
\begin{array}[c]{l}
\tilde f'''+\tilde f\tilde f''-\frac{2m}{m+1}\tilde f'^2=0, \\
\noalign{\vskip2mm}
\tilde f(0)=-a\sqrt{-\frac{m+1}2},~~
\tilde f'(0)=1,~~
\tilde f'(\infty)=0
\end{array}
\right.\label{pb-jde}
\end{equation}
and, from \cite{jde} (Theorem 1 and Proposition 1), we know that the problem $(\ref{pb-jde})$ admits exactly one concave solution, that this solution is bounded, and that if $\tilde \ell$ is the limit of $\tilde f$ at infinity, then $-a\sqrt{-\frac{m+1}2}<\tilde\ell<\sqrt{-\frac{m+1}2a^2+2}$. The proof is complete.\end{proof}

\medskip
\begin{remark} For $m<-1$, the uniqueness of the convex solution can be easily obtained by a direct way.
In fact, if $(\ref{equation})$-$(\ref{cond03})$ has a pair of distinct convex solutions $f_1$, $f_2$ and if $f''_1(0)>f''_2(0)$, the function $g=f_1-f_2$ satisfies $g(0)=0$, $g'(0)=g'(\infty)=0$ and $g''(0)>0$. It follows that $g'$ has a positive maximum at some point $s>0$ such that $g'(t)>0$ for $0<t<s$, but then $g(s)>0$ and we get
$$g'''(s)=f'''_1(s)-f'''_2(s)=-\frac{m+1}{2}f''_1(s)g(s)+m(f'_1(s)+f'_2(s))g'(s)>0,$$
which is a contradiction.
\end{remark}
\medskip

\begin{remark}\label{cvv}
For $m<-1$ and $a\leq 0$, if $f$ is a convex-concave solution of the problem $(\ref{equation})$-$(\ref{cond03})$, then 
\begin{equation*}
f''(0)>\max\left\{\frac{m+1}{2}a,\sqrt\frac{-2m}{3}\right\}.
\end{equation*}
 Indeed, if $t_1$ is the point such that $f'(t_1)=0$, we have $f''(t_1)>0$ and writing 
 equality $(\ref{i1})$ with $\alpha=0$ and $\beta=t_1$, we get $f''(0)>\frac{m+1}{2}a.$ 
Now, writing equality $(\ref{i2})$ with $\alpha=0$ and $\beta=\infty$ and taking into account the fact that $f<0$ on $(0,\infty)$ $($see {\rm Proposition} $\ref{prop<0})$, we obtain that $f''(0)>\sqrt\frac{-2m}{3}.$
\end{remark}

\medskip

\begin{lemma} \label{sep1} Let $m<-1$. As $s$ grows, the separatrix $S_0^-$ leaves the singular point $O$ to the right tangentially to $L$ and intersects successively the isoclines $Q_m(u,v)=0$, $P(u,v)=0$, the $u$-axis and the $v$-axis and remains decreasing  and under $L$. 

As $s$ grows, the separatrix $S_1^-$ leaves the singular point $O$ to the left tangentially to $L$ and remains decreasing and under $L$. \rm{(See figure $2$.)}
\end{lemma}
\begin{proof}
The behavior of $S_0^-$ is established in \cite{brighisari} for $u\geq 0$. To conclude, it is sufficient to remark that we have 
$$\frac{dv}{du}-\left(-\frac{m+1}{2}\right)=-\frac{u(3v+u)}{v-2u^2}>0$$ 
for $u,v<0$. Hence the straight line $L$ is a barrier for the vector field in the region $\{u<0\}\cap\{v<0\}$. \end{proof}
\medskip
\begin{center} 
\includegraphics[scale=0.38]{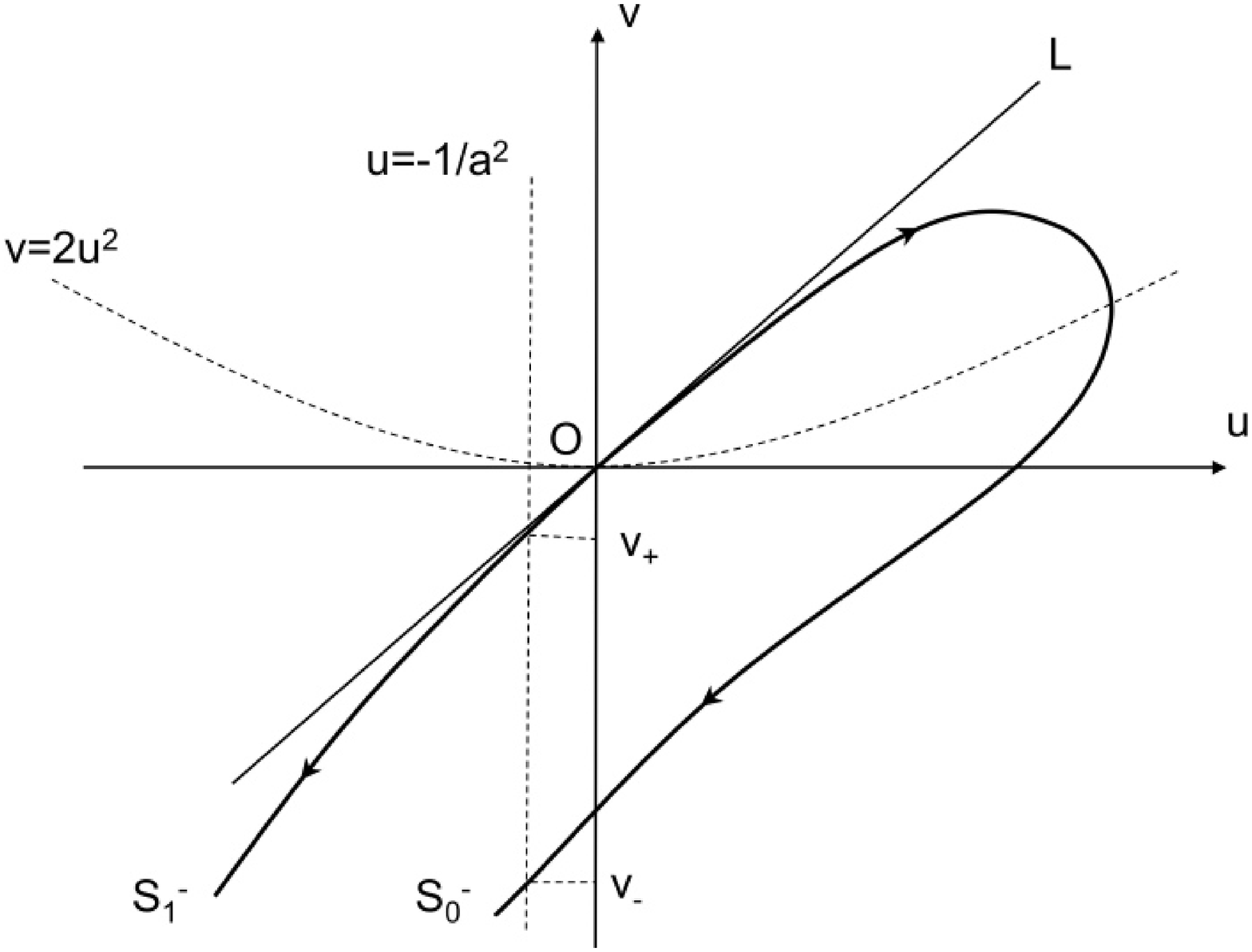}     
\end{center}

\medskip
\centerline{Figure 2 $:$ $m<-1$}
\medskip

\begin{theorem} Let $m<-1$. For every $a<0$, the problem $(\ref{equation})$-$(\ref{cond03})$ admits a unique convex-concave solution such that $\underset{t\rightarrow \infty}{\lim} f(t)=\ell<0$ and infinitely many convex-concave solutions such that $\underset{t\rightarrow \infty}{\lim} f(t)=0$. Moreover, all these solutions are negative and bounded.
\end{theorem}
\begin{proof} For $a<0$, let us denote by $f$ the solution of the initial value problem $({\mathcal P}_{m,a,b})$ and look at the corresponding trajectories $C_{a,b}$ of the plane system (\ref{system}) defined by (\ref{new_function}). From the phase portrait of (\ref{system}) described in Lemma \ref{sep1} (Fig 2) we have that each function $f$ corresponding to a trajectory that did not start from a point that is between the separatix $S_1^-$ and the separatrix $S_0^-$ for $u$ negative cannot be a solution of the problem  $(\ref{equation})$-$(\ref{cond03})$, because it vanishes (the trajectory goes to infinity) or becomes concave and decreasing (the trajectory goes through the domain $\{u<0\}\cap\{v>0\}$).
We also have that the straight line defined by $u=-\frac{1}{a^2}$ with $a<0$ intersects the separatrix $S_0^-$ at the point $\left(-\frac{1}{a^2},v_-\right)$ and the separatrix $S_1^-$ at the point $\left(-\frac{1}{a^2},v_+\right)$.

Consequently, for $b=v_-a^3$, $f$ is a solution of the problem $(\ref{equation})$-$(\ref{cond03})$ such that $f(t)\to\ell<0$ as $t\to\infty$, and for $b\in(v_+a^3,v_-a^3)$, $f$ is a solution of the problem $(\ref{equation})$-$(\ref{cond03})$ such that $f(t)\to 0$ as $t\to\infty$.
The details about the limits can be found in \cite{brighisari}.\end{proof}

\medskip

\begin{remark} Let us notice that the separatrix $S_1^-$ $($corresponding to $b=v_-a^3)$ gives, for $a<0$, the convex solution obtained in {\rm Theorem} $\ref{cv1}$.
\end{remark}
\subsection{The case $m=-1$.} 
Here the equation $(\ref{equation})$ reduces to $f'''=-f'^2$.  Now, if $f$ is a solution of the problem $(\ref{equation})$-$(\ref{cond03})$ then using the equality $(\ref{i2})$ with $\alpha=0$ and $\beta=\infty$, and Propositions $\ref{prop<0}$ and $\ref{f2}$, we obtain that $f''(0)=\sqrt{\frac{2}{3}}$. Hence the problem $(\ref{equation})$-$(\ref{cond03})$ has at most one solution. 

On the other hand, if $f$ is the solution of the problem $\left({\mathcal P}_{-1,a,\sqrt{\frac{2}{3}}}\right)$ on its right maximal interval of existence $[0,T)$ we get, using $(\ref{i2})$ with $\alpha=0$ and $\beta=t<T$,
$$\frac{1}{2}f''^2(t)+\frac{1}{3}f'^3(t)=0$$
implying that $f'<0$ on $[0,T)$. Therefore, we obtain
$$\forall t\in[0,T),~~~\frac{-f''(t)}{(-f'(t))^{\frac{3}{2}}}=-\sqrt{\frac{2}{3}}.$$
Integrating and using the fact that $f'(0)=-1$, we arrive to
$$\forall t\in[0,T),~~~\frac{1}{(-f'(t))^{\frac{1}{2}}}=\frac{t}{\sqrt{6}}+1$$
and $f'(t)=-\frac{6}{(t+\sqrt{6})^2}$. Integrating again, we finally deduce that, for $m=-1$, the problem $(\ref{equation})$-$(\ref{cond03})$ has exactly one solution given by
$$f(t)=a-\sqrt{6}+\frac{6}{t+\sqrt{6}}.$$
 
\subsection{The case $-1<m<0$.}
\begin{theorem} 
For $-1<m<0$ and for every $a\in\Bbb R$, the problem $(\ref{equation})$-$(\ref{cond03})$ admits a convex solution.
\end{theorem} 
\begin{proof} Let $b\geq 0$ and $f_b$ be the solution of the initial value problem $({\mathcal P}_{m,a,b})$. Denote by $[0,T_b)$ its right maximal interval of existence.
Let us remark first that $f_b$ exists as long as we have $f''_b>0$ and $f'_b<0$. Since $f''_b$  and $f'_b$ cannot vanish at the same point, it follows that there are only three possibilities

\vspace{3mm}
(a) $f''_ b$ becomes negative from a point such that $f'_b<0,$

(b) $f'_b$ becomes positive from a point such that $f''_b>0,$

(c) we always have $f'_b<0$ and $f''_b>0.$
\vspace{3mm}

\noindent As $f'_0(0)=-1<0$, $f''_0(0)=0$ and $f'''_0(0)=m<0$, we have that $f_0$ is of type (a), and by continuity it must be so for $f_b$ with $b>0$ small enough.

On the other hand, as long as $f''_b(t)\geq 0$ and $f'_b(t)\leq 0$, we have  $f'_b(t)\geq -1$ and $f_b(t)\geq -t+a$. Therefore
$(\ref{i1})$ leads to
\begin{align*}
f''_b(t) &= b-\frac{m+1}{2}(f_b(t)f'_b(t)+a)+\frac{3m+1}{2}\int_0^tf'^2_b(s)ds \cr
&\geq b-\frac{m+1}{2}(t+\vert a\vert+a)+\frac{3m+1}{2}\int_0^tf'^2_b(s)ds \cr
&\geq b-\frac{m+1}{2}(\vert a\vert+a)+C_mt 
\end{align*}
where $C_m=m$ if $-1<m\leq -\frac{1}{3}$, and $C_m=-\frac{m+1}{2}$ if $-\frac{1}{3}\leq  m<0$.
Integrating, we obtain
\begin{equation*}
0\geq f'_b(t) \geq -1+\left(b-\frac{m+1}{2}(\vert a\vert+a)\right)t+C_m\frac{t^2}{2}:=P_b(t).\label{p}
\end{equation*}
For $b$ large enough, the equation $P_b(t)=0$ has two positive roots $t_0<t_1$, and therefore, for such $a$ and $b$, we have that $f'_b(t_0)=0$ and $f''_b(t)>0$ for $t\leq t_0$, and $f_b$ is of type (b).

Defining $A=\left \{ b>0\, ; \, \text{$f_b$ is of type (a)}\right \}$ and
$B=\left \{ b>0\, ; \, \text{$f_b$ is of type (b)}\right \}$ we have that 
$A \neq \emptyset$, $B \neq \emptyset$ and $A \cap B = \emptyset$. Both $A$ and $B$ are open sets, so there exists a $b^*>0$ such that the corresponding solution $f_{*}$ of 
$\left({\mathcal P}_{m,a,b^*}\right)$ is of type (c) and is defined on the whole interval $[0,\infty)$. For this solution we have that $f'_{*}<0$ and $f''_{*}>0$ which implies that $f'_{*} \rightarrow l \leq 0$ as $t\rightarrow \infty$. If $l<0$ then $f_{*}(t)\sim lt$ as $t\to\infty$ and using $(\ref{i1})$ we easily obtain $f''_{*}(t)\sim ml^2t$ as $t\to\infty$, which contradicts the fact that $f''_{*}>0$. Then $l=0$ and $f_{*}$ is a convex solution of $(\ref{equation})$-$(\ref{cond03})$. This completes the proof.\end{proof}

\medskip
About the uniqueness of a convex solution, we only get it for $-1<m<-\frac{1}{3}$. To this aim, we will need the following Lemma. 

\medskip

\begin{lemma}\label{L^2}
Let $-1<m<-\frac{1}{3}$. If $f$ is solution of  the problem $(\ref{equation})$-$(\ref{cond03})$, then we have $f(t)f'(t)\to 0$ as $t\to\infty$ and 
\begin{equation}
f''(0)=\frac{m+1}{2}a-\frac{3m+1}{2}\int_0^\infty f'^2(s)ds.\label{L2}
\end{equation}
\end{lemma}
\begin{proof} If $f$ is bounded, it clearly holds that $f(t)f'(t)\to 0$ as $t\to\infty$, and $(\ref{L2})$ follows immediatly from $(\ref{f''})$ and the equality $(\ref{i1})$ written with $\alpha=0$ and $\beta=\infty$.

If $f$ is unbounded, then, as $t\to\infty$, either $f$ increases to $\infty$ or decreases to 
$-\infty$, and in both cases we have $f(t)f'(t)>0$ for $t$ large enough. Now, 
using the equality $(\ref{i1})$ with $\alpha=0$ and $\beta=t>0$, for large $t$, we get
\begin{equation}
0<-\frac{3m+1}{2}\int_0^tf'^2(s)ds\leq f''(0)-f''(t)-\frac{m+1}{2}a.
\label{L2'}
\end{equation}
Thus $(\ref{f''})$ implies that the integral in the relation $(\ref{L2'})$ has a limit as $t\to\infty$. Coming back to $(\ref{i1})$, we get that $f(t)f'(t)\rightarrow l\geq 0$ as $t\rightarrow \infty$. If $l>0$, then $f^2(t)\sim 2lt$ as $t\to\infty$, and hence
$$f'^2(t)\sim\frac{l^2}{f^2(t)}\sim\frac{l}{2t}~~\hbox{ as }~t\rightarrow \infty$$
which contradicts the fact that the integral of $f'^2$ over $[0,\infty)$ is finite. Therefore $l=0$, and using again $(\ref{i1})$  we obtain $(\ref{L2})$.\end{proof}

\medskip
\begin{proposition} For $-1<m<-\frac{1}{3}$ and $a\in\Bbb R$, the problem $(\ref{equation})$-$(\ref{cond03})$ admits at most one convex solution.
\end{proposition}
\begin{proof}
Let  us suppose that $f$ is a  
convex solution of (\ref{equation})-(\ref{cond03}). By Proposition \ref{prop<0}, such a solution is decreasing, and we can define the function $v=v(y)$ by
$$\forall t\geq 0,\quad v(f(t))=f'(t).$$
If $\ell\in[-\infty,a)$ is the limit of $f(t)$ as $t \to \infty$, then $v$ is defined on $(\ell,a]$, is negative and we have
\begin{align*}
f''(t)&=v(f(t))v'(f(t)), \\
f'''(t)&=v(f(t))v'^2(f(t))+v^2(f(t))v''(f(t)).
\end{align*}
Then, (\ref{equation}) leads to
\begin{equation}
\forall y\in (\ell,a],\quad v''=-\frac{1}{v}\left (v'+\frac{m+1}{2}y\right)v'+m
\label{eq-v}
\end{equation}
and we have
$$v(\ell):=\underset{y\rightarrow \ell}{\lim}\, v(y)
=\underset{t\rightarrow \infty}{\lim}\, f'(t)=0, \quad v(a)=v(f(0))=f'(0)=-1.$$
Let us now suppose that there are two convex solutions $f_{1}$ and $f_{2}$ of the problem (\ref{equation})-(\ref{cond03}) and let $\ell_{i}\in[-\infty,a)$ be the limit of $f_{i}$ at infinity for $i=1,2.$ We obtain two solutions of equation (\ref{eq-v}), $v_{1}$ and $v_{2}$  defined respectively on $(\ell_{1},a]$ and $(\ell_{2},a]$ such that
$$v_{1}(\ell_{1})=v_{2}(\ell_{2})=0 \quad \text{and}\quad v_{1}(a)=v_{2}(a)=-1.$$
We will suppose now that $\ell_2\leq\ell_1$ and prove that $v_2\leq v_1$ on $(\ell_1,a]$.
If there exists a point $y$ in $(\ell_1,a]$ such that $v_1(y)<v_2(y)$, then as $v_1(\ell_1)-v_2(\ell_1)\geq 0$ and $v_1(a)-v_2(a)=0$, the function $v_1-v_2$ admits a negative minimum at some point $x$ in $(\ell_1,a)$. For this point $x$ we have that $v_{1}(x)<v_{2}(x)$, $v'_{1}(x)=v'_{2}(x)$ and $v''_{1}(x)\geq v''_{2}(x)$. We also have
\begin{equation}
v''_{1}(x)-v''_{2}(x)=\left (\frac{1}{v_{2}(x)}-\frac{1}{v_{1}(x)} \right )
\left (v'_{1}(x)+\frac{m+1}{2}x\right)v'_{1}(x) 
\label{eq-v1}
\end{equation}
and
\begin{equation}
\left (v'_{1}(x)+\frac{m+1}{2}x\right)v'_{1}(x)=
\left (f''_{1}(s)+\frac{m+1}{2}f_{1}(s)f'_{1}(s)\right)\frac{f''_{1}(s)}{f'^2_{1}(s)} 
\label{eq-v2}
\end{equation}
with $s$ such that $x=f_{1}(s)$. Thanks to Lemma \ref{L^2}, we can write (\ref{i1}) with $\alpha=s$ and $\beta=\infty$ to get
\begin{equation*}
f''_{1}(s)+\frac{m+1}{2}f_{1}(s)f'_{1}(s)=-\frac{3m+1}{2}\int_{s}^{\infty}f'^2_{1}(t)dt>0. \end{equation*}
Using this inequality and the fact that $f''_{1}(s)>0$, we deduce from (\ref{eq-v1}) and (\ref{eq-v2}) that 
$v''_{1}(x)<v''_{2}(x)$ and obtain a contradiction. Therefore we have $v_2\leq v_1$ on $(\ell_1,a]$ 
and
$$\int_{0}^{\infty}f'^2_{2}(t)dt=\int^{a}_{\ell_2}(-v_{2}(y))dy\geq\int^{a}_{\ell_1}(-v_{2}(y))dy
\geq \int^{a}_{\ell_1}(-v_{1}(y))dy=\int_{0}^{\infty}f'^2_{1}(t)dt.$$
From Lemma \ref{L^2} we have
$$f''_i(0)=\frac{m+1}{2}a-\frac{3m+1}{2}\int_{0}^{\infty}f'^2_{i}(t)dt$$ and thus $f''_1(0)\leq f''_2(0)$. If $f''_1(0)<f''_2(0)$ then $v_1(a)v'_1(a)<v_2(a)v'_2(a)$ and, as $v_1(a)=v_2(a)=-1$, this leads to $v'_1(a)>v'_2(a)$ that is a contradiction with the fact that $v_1\geq v_2$ on $(\ell_1,a]$. Hence $f''_1(0)=f''_2(0)$ and $f_1=f_2$.\end{proof}

\medskip

\begin{proposition}
Let $-1<m\leq -\frac{1}{2}$. If $a\leq 0$, then the problem $(\ref{equation})$-$(\ref{cond03})$ has no convex-concave solution.
\end{proposition}
\begin{proof} 
Let us suppose that $f$ is a convex-concave solution of $(\ref{equation})$-$(\ref{cond03})$.
By Proposition  \ref{prop<0} we have that $f>0$ at infinity, and thus there exists $\tau>0$ such that $f(\tau)=0$ and $f'(\tau)>0$. Using $(\ref{i3})$ written with $\alpha=\tau$ and $\beta=t_n$ where $(t_n)$ is the sequence defined in Proposition \ref{f2}, we get 
\begin{align*}
0&\geq\lim_{n\to \infty}(2m+1)\int_\tau^{t_n}f'^2(t)f(t)dt\\
&=\lim_{n\to \infty}\left(\frac12f'^2(\tau)+\frac{m+1}{2}f'(t_n)f^2(t_n)\right)\geq\frac12f'^2(\tau)>0
\end{align*}
and a contradiction.\end{proof}

\medskip 

\begin{remark}
For $m=-\frac{1}{3}$, let $f$ be the solution of $({\mathcal P}_{m,a,b})$ on its right maximal interval of existence $[0,T)$. Integrating twice leads to the Riccati equation
\begin{equation}
f'+\frac{1}{6}f'^2=\left(b-\frac{a}{3}\right)t-1+\frac{a^2}{6}. \label{rc}
\end{equation}
Choosing $b=\frac{a}{3}$, equation $(\ref{rc})$ has explicit particular solutions, and can be solved. We then obtain that 

\begin{equation}
f(t)=\sqrt{a^2-6}\cdot\frac{(a+\sqrt{a^2-6})e^{\frac t6\sqrt{a^2-6}}+(a-\sqrt{a^2-6})e^{-\frac t6\sqrt{a^2-6}}}{(a+\sqrt{a^2-6})e^{\frac t6\sqrt{a^2-6}}-(a-\sqrt{a^2-6})e^{-\frac t6\sqrt{a^2-6}}}\label{rc1}
\end{equation}
when $a^2>6$,
$$f(t)=\sqrt{6-a^2}\cdot\rm{cotan}\left(\textstyle\frac t6\sqrt{6-a^2}+\rm{arccotan}\frac{a}{\sqrt{6-a^2}}\right)$$
when $a^2<6$ and 
\begin{equation}
f(t)=\frac{6}{t+a}\label{rc2}
\end{equation}
when $a^2=6$. It is then easy to see that $T$ is finite if $a<\sqrt{6}$, that $T=\infty$  if  $a\geq\sqrt{6}$, and that, in this latter case, the function $f$ given by  $(\ref{rc1})$ for $a>\sqrt{6}$ and by $(\ref{rc2})$ for $a=\sqrt{6}$ is a solution of $(\ref{equation})$-$(\ref{cond03})$.
\end{remark}

\medskip

In the remainder of this section we will concentrate our efforts to the case $m>0$.
In this case, we know from Proposition \ref{prop>0} that we must have $a>0$ and that the solutions of the problem $(\ref{equation})$-$(\ref{cond03})$ cannot vanish, this allows us to consider the dynamical system (\ref{system}). We will consider successively the cases $0<m<1$, $m=1$ and $m>1$, but before let us remark that we can improve slightly the result of Proposition \ref{prop>0}. Actually, if $m>0$ then, for $a\leq\frac{2}{\sqrt{m+1}}$, the problem $(\ref{equation})$-$(\ref{cond03})$ has no solutions. Indeed, if $f$ is a solution of $(\ref{equation})$-$(\ref{cond03})$, then $f$ is bounded and if we write $(\ref{i1})$ with $\alpha=t$ and $\beta=\infty$, we get
$$-f''(t)-\frac{m+1}{2}f(t)f'(t)=\frac{3m+1}{2}\int_t^\infty f'^2(s)ds>0.$$
Integrating we obtain $$-f'(t)-1-\frac{m+1}{4}(f^2(t)-a^2)>0$$ and by taking the limit as $t\to\infty$, we get $-1-\frac{m+1}{4}(\ell^2-a^2)>0$,
where $\ell$ is the limit of $f$ at infinity. This gives $a>\frac{2}{\sqrt{m+1}}$, and also $0\leq \ell<\sqrt{a^2-\frac{4}{m+1}}$.
\medskip

\subsection{The case $0<m<1$.}
\begin{theorem} 
Let $0<m<1$, then there exists $0<a_1^*<a_2^*$ such that the problem $(\ref{equation})$-$(\ref{cond03})$ admits no solutions for $0<a<a_1^*$, one unique solution for $a=a_1^*$, multiple solutions for $a_1^*<a<a_2^*$, two solutions for $a=a_2^*$ and a unique solution for $a_2^*<a$. Moreover, all these solutions are convex, decreasing and positive.
\end{theorem}
\begin{proof} In this case, the point $A$ is an unstable focus. The separatrix $S_1^+$ leaves the point $A$ or a limit cycle surrounding $A$ that stays in the domain $\{u<0\}\cap\{v>0\}$, turning clockwise and cut the isoclinic curve $P(u,v)=0$ for its last time, then the isoclinic curve $Q_m(u,v)=0$ and goes to the point $O$ (Fig 3). For the proof, refer to \cite{brighisari}. 

For $a>0$, we denote by $f$ the solution of the initial value problem $({\mathcal P}_{m,a,b})$ and look at the corresponding trajectories $C_{a,b}$ of the plane system (\ref{system}) defined by (\ref{new_function}).

Let us now consider the straight line $u=-\frac{1}{a^2}$ for $a>0$ and let $u_1^*$ (i.e. $a_1^*=\sqrt{-\frac{1}{u_1^*}} $) be the abscissa of the point at which the separatrix $S_1^+$ crosses the isoclinic curve $P(u,v)=0$ for its last time and $u_2^*$ (i.e. $a_2^*=\sqrt{-\frac{1}{u_2^*}}$) be the abscissa of the point at which the separatrix $S_1^+$ crosses the isoclinic curve $P(u,v)=0$ for its penultimate time.
\medskip
\begin{center}
\includegraphics[scale=0.55]{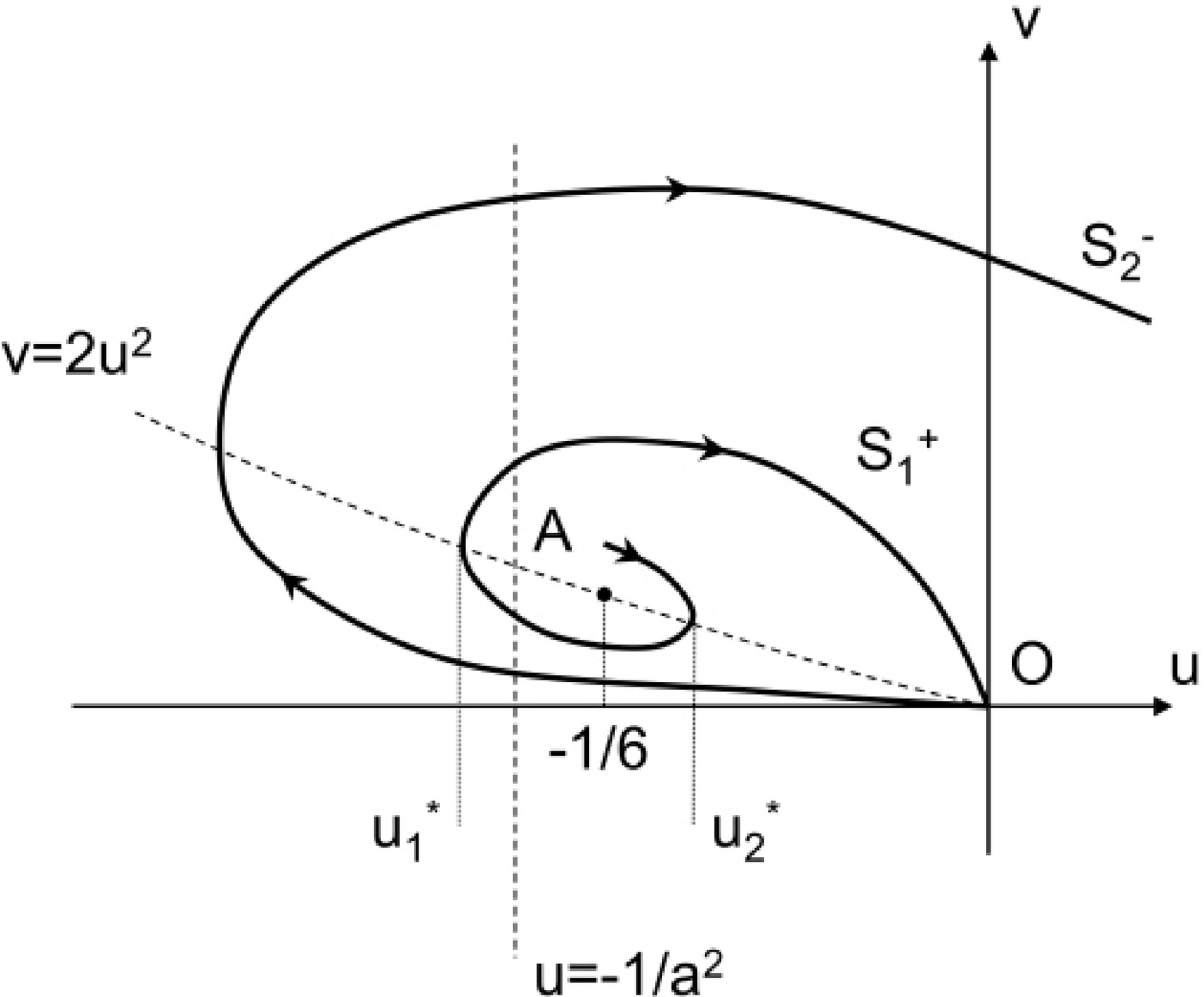}     
\end{center}

\medskip
\centerline{Figure 3 $:$ $0<m<1$}
\medskip
Looking at the phase portrait of (\ref{system}) we see immediatly that if $u<u_1^*$, then $f$ is not a solution of the problem (\ref{equation})-(\ref{cond03}) because the corresponding trajectories crosses the $v$-axis (meaning that $f$ becomes increasing) or going to infinity (meaning that $f$ vanishes). 

Moreover, for the same purpose, the function $f$ corresponding to parts of the separatrix $S_1^+$, to parts of the limit cycles surrounding $A$ or to parts of the trajectories inside the limit cycles, is a convex solution of (\ref{equation})-(\ref{cond03}), because $u$ and $v$ remain  positive.

Consequently, for $u=u_1^*$ there is only one solution, for $u_1^*<u<u_2^*$ there are multiple solutions, for $u=u_2^*$ there are two solutions,  and for $u_2^*<u<0$ there is only one solution.\end{proof}

\subsection{The case $m=1$.}
The equation $(\ref{equation})$ is
$$f'''+ff''-f'^2=0.$$
If we set $f=g+k$ with $k\in\Bbb R$, then $f$ is a solution of $(\ref{equation})$-$(\ref{cond03})$ if and only if $g$ satisfies
\begin{equation}
g'''+kg''=g'^2-gg''\label{g-eq}
\end{equation}
and $g(0)=a-k$, $g'(0)=-1$, $g'(\infty)=0$.
Looking for functions $g$ such that both hand sides of $(\ref{g-eq})$ vanish, we get that, for $a\geq 2$, the functions $f_1,f_2:[0,\infty)\to\Bbb R$ given by
$$f_i(t)=k_i+\frac{1}{k_i}e^{-k_i t}~~~\hbox{ for }~i=1,2$$
with $k_1=\frac12(a-\sqrt{a^2-4})$ and $k_2=\frac12(a+\sqrt{a^2-4})$ are convex solutions of $(\ref{equation})$-$(\ref{cond03})$. 
\medskip

\subsection{The case $m>1$.}
\begin{theorem} 
Let $m>1$, then there exists $0<a_1^*<a_2^*$ such that the problem $(\ref{equation})$-$(\ref{cond03})$ admits no solutions for $0<a<a_1^*$, one unique solution for $a=a_1^*$ that is convex and such that $\underset{t\rightarrow \infty}{\lim} f(t)=\ell>0$, two convex solutions that verifies $\underset{t\rightarrow \infty}{\lim} f(t)=\ell>0$ and infinitely many convex solutions such that $\underset{t\rightarrow \infty}{\lim} f(t)=0$ for $a_1^*<a\leq a_2^*$, one concave-convex solutions and one convex solution such that $\underset{t\rightarrow \infty}{\lim} f(t)=\ell>0$ and infinitely many concave-convex or convex solutions with $\underset{t\rightarrow \infty}{\lim} f(t)=0$ for $a_2^*<a$. 
All these solutions are decreasing and positive.
\end{theorem}
\begin{proof}
Let us recall that the point $A$ is an unstable focus for $1<m<\frac{4}{3}$, a stable focus for $\frac{4}{3}<m\leq \frac{4+2\sqrt{6}}{3}$ or a stable node for $m \geq \frac{4+2\sqrt{6}}{3}$. 
Moreover, for $1<m<\frac{4}{3}$, as $A$ is unstable there exists at least one cycle surrounding it. The behavior of the separatrices $S_1^+$ and $S_2^-$ is established in \cite{brighisari} and described by the figure 4. 

For $a>0$, let us  denote again by $f$ the solution of the initial value problem $({\mathcal P}_{m,a,b})$ and look at the corresponding trajectories $C_{a,b}$ of the plane system (\ref{system}) defined by (\ref{new_function}). 
 
Consider now the straight line $u=-\frac{1}{a^2}$ for $a>0$ and let $u_1^*$  (i.e. $a_1^*=\sqrt{-\frac{1}{u_1^*}} $) be the abscissa of the point at which the separatrix $S_1^+$ crosses the isoclinic curve $P(u,v)=0$  and $u_2^*$  (i.e. $a_2^*=\sqrt{-\frac{1}{u_2^*}} $) be the abscissa of the point at which the separatrix $S_1^+$ crosses the $u$-axis.

Let $\mathcal D$ be the bounded domain delimited by the $v$-axis and the part of the separatrix $S_1^+$ included in $\{u<0\}$.

The phase portrait of (\ref{system}) gives us immediatly that every trajectory that starts from a point with negative abscissa and outside of $\mathcal D$ crosses the $v$-axis (meaning that $f$ becomes increasing) or is going to infinity (meaning that $f$ vanishes). Then, the corresponding function $f$ is not a solution of the problem (\ref{equation})-(\ref{cond03}). 
\begin{center}
\includegraphics[scale=0.45]{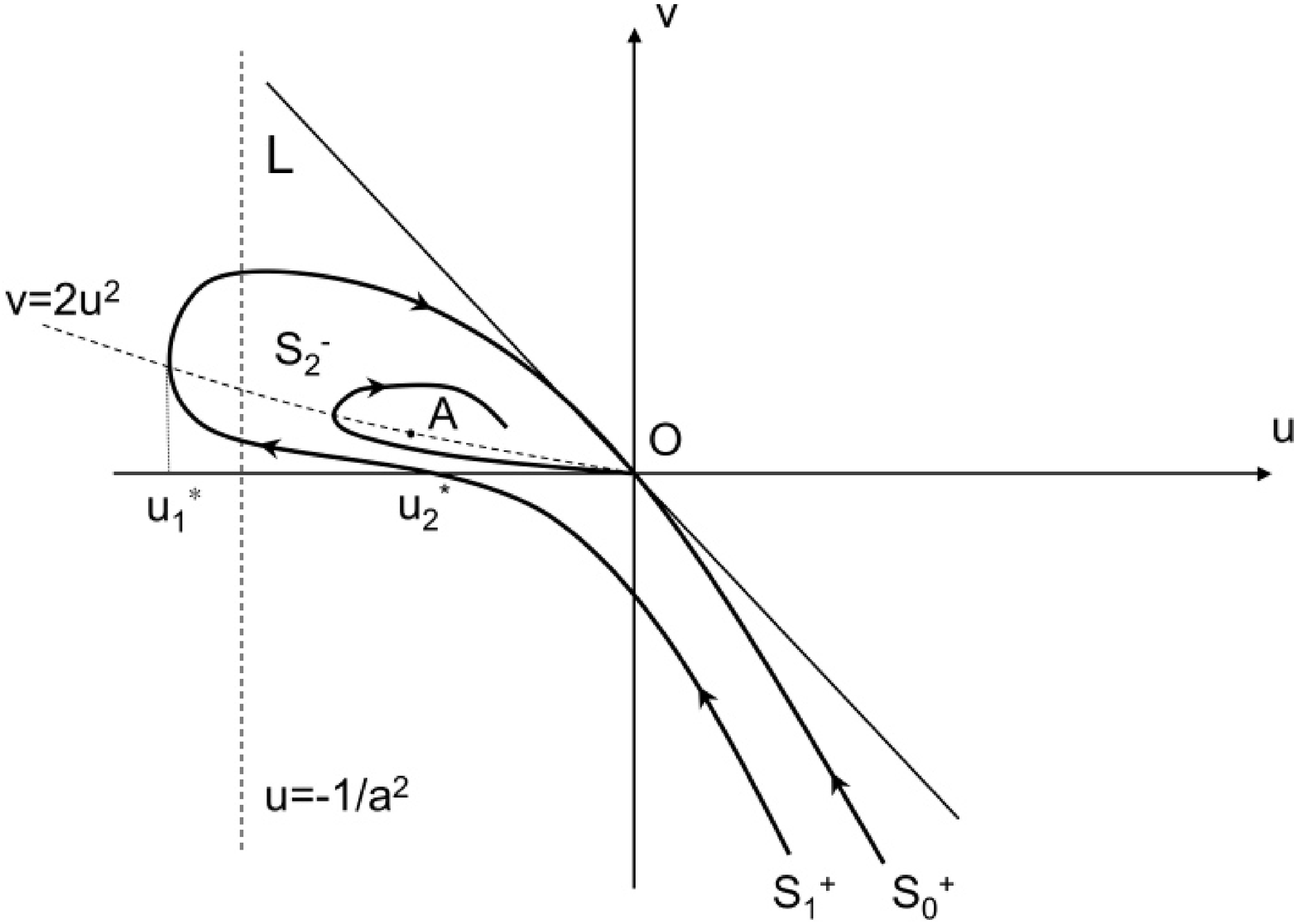}     
\end{center}

\medskip
\centerline{Figure 4 $:$ $m>1$}
\medskip

As shown in \cite{brighisari}, using the Poincar\'e-Bendixson Theorem, every trajectory that enters in $\mathcal D$ must have the point $O$, the point $A$ or a limit cycle surrounding $A$ that is not crossing the $u$-axis for $\omega$-limit set. This means that every function $f$ that corresponds to a phase curve that starts from a point in $\bar{ \mathcal D}$ is a convex solution of (\ref{equation})-(\ref{cond03}), because $u$ and $v$ remain  positive.
 This leads to the following behavior.

For $u<u_1^*$, there is no solution of $(\ref{equation})$-$(\ref{cond03})$.

For $u=u_1^*$, there is a unique solution of $(\ref{equation})$-$(\ref{cond03})$ that is given by a part of the separatrix $S_1^+$. This solution is convex, decreasing and such that $f(t)\to\ell>0$ as $t\to\infty$.

For $u_1^*<u\leq u_2^*$, there are two solutions given by parts of the separatrix $S_1^+$ that are convex and such that $f(t)\to\ell>0$ as $t\to\infty$, and infinitely many convex solutions given by the trajectories starting from a point in $\mathcal{D}$ with an abscissa equal to $-\frac{1}{a^2}$ and such that $f(t)\to 0$ as $t\to\infty$.

For $u_2^*<u<0$, there is one concave-convex solution and one convex solution given by parts of the separatrix $S_1^+$ that both verify $f(t)\to\ell>0$ as $t\to\infty$, and infinitely many concave-convex or convex solutions such that $f(t)\to 0$ as $t\to\infty$ corresponding to the trajectories starting from a point in $\mathcal{D}$ with an abscissa equal to $-\frac{1}{a^2}$.

For the proofs of the limits as $t$ goes to infinity, refer to \cite{brighisari}.\end{proof}

\section{Conclusion}

In this paper we have studied a boundary value problem involving a third
order ordinary differential equation. The solutions of this problem are
similarity solutions of problems related to the phenomenon of high frequency
excitation of liquid metal systems in an antisymmetric magnetic field, within the framework of boundary layer approximation.

This study can be compared to the one made in \cite{brighisari} which
consider the same equation but where the solutions are supposed to start
increasing from $0$ instead of decreasing.

We have established several results in both cases $m<0$ and $m>0$.
Nevertheless, some interesting open questions still remain
\begin{itemize}
\item for $m<-1$ and $a\geq 0$, is there convex-concave solutions?
\item for $-\frac{1}{3}<m<0$ and $a\in \mathbb R$, is the convex solution unique?
\item for $-1<m\leq-\frac{1}{2}$ and $a>0$, is there convex-concave solutions?
\item for $-\frac{1}{2}<m<0$ and $a\in \mathbb R$, is there convex-concave solutions?
\end{itemize}

\section*{Acknowledgement} 
The first author wish to thanks A. Fruchard and T. Sari for stimulating discussions. 

The second author would like to thank M. Guedda and Z. Hammouch for enjoyable discussions and for bringing the subject of this paper to his attention as well as G. Karch and all the organizers for the wonderful time he spent in Bedlewo during the conference.


\begin{thebibliography}{99}   

\bibitem{ake} S. Akesbi, B. Brighi \& J.-D. Hoernel, Steady free convection in a bounded and saturated porous medium, Proceedings of the Swiss-Japanese Seminar on
Elliptic and Parabolic Issues in Applied Sciences, Z\"urich, December 2004,
pp. 1-17. World Scientific Publishing Co. Pte. Ltd., 2006.

\bibitem{aly} E. H. Aly, L. Elliott \& D. B. Ingham, Mixed convection boundary-layer flows over a vertical surface embedded in a porous medium, Eur. J. Mech. B Fluids 22 (2003), pp. 529-543.

\bibitem{banks1}W. H. H. Banks, Similarity solutions of the boundary layer equations
for a stretching wall, J. de M\'echan. Th\'eor. et  Appl. 2 (1983), pp. 375-392.                                                  

\bibitem{brighi03}Z. Belhachmi, B. Brighi \& K. Taous, On the concave
solutions of the Blasius equation, Acta Math. Univ. Comenianae, Vol. LXIX, 2 (2000), 
pp. 199-214.

\bibitem{brighicr}Z. Belhachmi, B. Brighi \& K. Taous, Solutions similaires pour un 
probl\`eme de couche limite en milieux poreux, C. R. M\'ecanique 328 (2000), 
pp. 407-410.

\bibitem {brighi02}Z. Belhachmi, B. Brighi \& K. Taous, On a family of
differential equations for boundary layer approximations in porous media,
Euro. Jnl of Applied Mathematics, Vol. 12, 4 (2001),
pp. 513-528.

\bibitem {brighi04}Z. Belhachmi, B. Brighi, J. M. Sac-Ep\'ee \& K. Taous,
Numerical simulations of free convection about a vertical flat plate embedded
in a porous medium, Computational Geosciences, vol. 7 (2003), pp. 137-166.

\bibitem{bla}H. Blasius, Grenzschichten in Fl\"ussigkeiten mit kleiner Reibung,
Z. Math. Phys. 56 (1908), pp. 1-37.

\bibitem {brighi01}B. Brighi, On a similarity boundary layer equation,
Zeitschrift f\"ur Analysis und ihre Anwendungen, vol. 21, 4 (2002), pp. 931-948.

\bibitem {bfs}B. Brighi, A. Fruchard \& T. Sari, On the Blasius problem. Preprint.

\bibitem {heat_flux}B. Brighi, J.-D. Hoernel, On similarity solutions for boundary layer flows with prescribed heat flux. Mathematical Methods in the Applied Sciences, vol. 28, 4 (2005) pp. 479-503.

\bibitem{gaeta} B. Brighi, J.-D. Hoernel, Recent advances on similarity solutions arising during free convection, Progress in Nonlinear Differential Equations and Their Applications, Vol. 63, pp. 83-92, Birkh\"auser Verlag Basel/Switzerland, 2005.

\bibitem{equiv}B. Brighi, J.-D. Hoernel, Asymptotic behavior of the unbounded 
solutions of some boundary layer equation, Archiv der Mathematik, Vol. 85, 2 (2005), pp. 161-166.

\bibitem {aml}B. Brighi, J.-D. Hoernel, On the concave and convex solutions of mixed convection boundary layer approximation in a porous medium, Applied Mathematics Letters, Vol 19, 1 (2006), pp. 69-74.

\bibitem {jde}B. Brighi, J.-D. Hoernel, On a general similarity boundary layer equation. Preprint.

\bibitem {brighisari}B. Brighi, T. Sari, Blowing-up coordinates for a
similarity boundary layer equation. Discrete and Continuous Dynamical Systems 
(Serie A), Vol. 12, 5 (2005), pp. 929-948.

\bibitem {pop1}M. A. Chaudhary, J.H. Merkin \& I. Pop, Similarity solutions in
free convection boundary-layer flows adjacent to vertical permeable surfaces
in porous media: I prescribed surface temperature, Eur. J. Mech. B-Fluids, 14 (1995), pp. 217-237.

\bibitem {pop}M. A. Chaudhary, J.H. Merkin \& I. Pop, Similarity solutions in
free convection boundary-layer flows adjacent to vertical permeable surfaces
in porous media: II prescribed surface heat flux, Heat and Mass Transfer 30,
Springer-Verlag (1995), pp. 341-347.

\bibitem{cheng}P. Cheng, W. J. Minkowycz, Free-convection about a vertical flat
plate embedded in a porous medium with application to heat transfer from a dike,
J. Geophys. Res. 82 (14) (1977), pp. 2040-2044.

\bibitem{coppel} W. A. Coppel, On a differential equation of boundary layer theory, Phil. Trans. Roy. Soc. London, Ser A 253, pp. 101-136 (1960).

\bibitem{crane}L. E. Crane, Flow past a stretching plane, Z. Angew. Math. Phys.
21 (1970), pp. 645-647.


\bibitem{falk}V. M. Falkner, S. W. Skan, Solutions of the boundary layer equations, Phil.
Mag., 7/12 (1931), pp. 865-896.


\bibitem{guedda} M. Guedda, Nonuniqueness of solutions to differential equations for boundary layer approximations in porous media, C. R. M\'ecanique, 330 (2002), pp. 279-283.

\bibitem{guedda1}M. Guedda, Similarity solutions of differential equations for boundary
layer approximations in porous media, Z. angew. Math. Phys. 56 (2005), pp. 749-762.

\bibitem{guedda2} M. Guedda, Multiple solutions of mixed convection boundary-layer approximations in a porous medium, Applied Mathematics Letters, Vol 19, 1 (2006), pp. 63-68.

\bibitem{gup}P. S. Gupta, A. S. Gupta, Heat an mass transfer on a stretching sheet
with suction or blowing, Can. J. Chem. Eng. 55 (1977), pp. 744-746.

\bibitem {hart}P. Hartmann, Ordinary Differential Equations. Wiley, New-York (1964).



\bibitem{ing}D. B. Ingham, S. N. Brown, Flow past a suddenly heated vertical plate
in a porous medium, J. Proc. R. Soc. Lond. A 403 (1986), pp. 51-80.

\bibitem{ish1}N. Ishimura, T. K. Ushijima, An elementary approach to the analysis of exact solutions for the Navier-Stokes stagnation flows with slips, Arch. Math. 82 (2004), pp. 432-441.




\bibitem{mag} E. Magyari, B. Keller, Exact solutions for self-similar boundary-layer flows induced by permeable stretching wall. Eur. J. Mech. B-Fluids 19 (2000), pp. 109-122.



\bibitem{mof}H. K. Moffatt, High-frequency excitation of liquid metal systems, IUTAM Symposium: Metallurgical Application of Magnetohydrodynamics, (1982) Cambridge.

\bibitem{nazar} R. Nazar, N. Amin \& I. Pop, Unsteady mixed convection boundary-layer flow near the stagnation point on a vertical surface in a porous medium, Int. J. Heat Mass Transfer  47 (2004), pp. 2681-2688.






\bibitem{stu} J. T. Stuart, Double boundary layers in oscillatory viscous flow, J. Fluid. Mech. 24 (1966), pp. 673-687. 



\bibitem{utz}W. R. Utz, Existence of solutions of a generalized Blasius equation, J.
Math. Anal. Appl. 66 (1978), pp. 55-59.


\bibitem{wang}J. Wang, W. Gao \& Z. Zhang, Singular nonlinear boundary value problems arising in boundary layer theory, J. Math. Anal. Appl. 233 (1999), pp. 246-256.



\bibitem{wood}R. A. Wooding, Convection in a saturated porous medium at large Rayleigh number or Peclet number, J. Fluid. Mech., 15 (1963), pp. 527-544.

\bibitem{yang} G.C. Yang, Existence of solutions to the third-order nonlinear differential
equations arising in boundary layer theory, Appl. Math. Lett. 16 (6) (2003),
pp. 827-832. 

\bibitem{yang2} G.C. Yang, A note on $f'''+ff''+\lambda(1 - f'^2)=0$ with $\lambda \in (-1/2, 0)$
arising in boundary layer theory, Appl. Math. Lett. 17 (11) (2004), pp. 1261-1265. 

\end{thebibliography}
\end{document}